\documentclass[preprint,12pt]{elsarticle}

\usepackage{amssymb}
\usepackage{amscd,amsmath,amsfonts,graphicx}
\usepackage{amsthm}

\newcommand{\R}{\mathbb{R}}
\newcommand{\C}{\mathbb{C}}

\newcommand{\N}{\mathbb {N}}
\newcommand{\E}{\mathsf{E}}

\begin{document}
\baselineskip=22pt \centerline{\large \bf Integral Representations
of Three Novel Multiple Zeta  Functions} \centerline{\large \bf
 for Barnes Type: A Probabilistic Approach}

\vspace{1cm} \centerline{Gwo Dong Lin $^a$ and Chin-Yuan Hu $^b$}
\centerline{$^a$ Academia Sinica, Taiwan}
\centerline{$^b$ National
Changhua University of Education, Taiwan}

 \vspace{1cm} \noindent {\bf Abstract.} Integral representation is one of the  powerful tools for studying
 analytic continuation of the zeta functions. It is known that  Hurwitz zeta
function generalizes the famous Riemann zeta function which plays an
important role in analytic number theory. They both have several
multiple versions in the literature.
 In this paper, we introduce three novel multiple
 zeta functions for Barnes type and study their integral representations
through hyperbolic probability distributions given by Pitman and Yor
({Canad. J. Math.}, {55} (2003), 292--330). The analytically
continued properties of the three multiple zeta functions are also
investigated. The first one generalizes the well-known Barnes
multiple zeta function, while the second and the third, unlike the
previous results, can extend analytically to entire functions in the
whole complex plane. These are somewhat surprising results.
\vspace{0.3cm}\\
\hrule
{\bf Running title:} Three novel Barnes
 multiple zeta
 functions\\
{\bf Postal addresses:}\\  Gwo Dong Lin, Institute of Statistical
Science, Academia Sinica,  No.\,128, Section 2, Academia Road,
Nankang District, Taipei 11529, Taiwan (ROC). (E-mail:
gdlin@stat.sinica.edu.tw)\\
Chin-Yuan Hu,  National Changhua
University of Education,  No.\,1, Jinde Road, Changhua 50007, Taiwan
(ROC).
(E-mail: buhuua@gmail.com)\\

\newpage
\noindent{\bf 1. Introduction}
\newcommand{\bin}[2]{
   \left(
     \begin{array}{@{}c@{}}
         #1  \\  #2
     \end{array}
   \right)          }

Consider the Hurwitz zeta function
\begin{eqnarray}\zeta(s;b)=\sum_{n=0}^{\infty}\frac{1}{(n+b)^s},
\end{eqnarray}
where the parameter $b>0$ and the complex variable $s\in\C$ has the
real part ${\rm Re}(s)>1.$ When $b=1$, $\zeta(s;b)$ reduces to the
famous Riemann zeta function
\begin{eqnarray}\zeta(s)\equiv\zeta(s;1)=\sum_{n=1}^{\infty}\frac{1}{n^s},\quad  {\rm Re}(s)>1,
\end{eqnarray}
which plays an important role in analytic number theory (see, e.g., Apostal
\cite{Apo1976}).
 It is known that each of both zeta functions in (1) and (2) can
extend to a meromorphic function in the whole complex plane $\C$ and
has a simple pole at $s=1$ where it has residue 1.

As for multidimensional generalizations of (1) and (2), there are
various multiple ($r$-tuple or $r$-fold) zeta functions proposed in
the literature (see, e.g., the survey paper by Matsumoto
\cite{Mat2010}, the thesis by Onozuka \cite{Ono2014},  and the
references therein), where $r\ge 1$
is an integer. We recall a few of the results in the following: \\
(i) Euler--Zagier multiple (Riemann) zeta function (\cite{Eul1776,
Hof1992, Mat2002a, MSV2016, Zag1994, Zha2000})
\[\zeta(s_1,s_2,\ldots,s_r)=\sum_{0<n_1<n_2<\cdots<n_r}\frac{1}{n_1^{s_1}n_2^{s_2}\cdots n_r^{s_r}},\quad \sum_{k=j}^r{\rm Re}(s_k)>r-j+1,\ \ 1\le j\le r;\]
(ii) Multiple Hurwitz  zeta function (\cite{Kam2006, KM2008, MV2017,
MS2006})
\[\zeta(s_1,s_2,\ldots,s_r;b)=\sum_{0\le n_1<n_2<\cdots<n_r}\prod_{j=1}^r\frac{1}{(n_j+b)^{s_j}}, \quad b> 0,\quad \sum_{k=j}^r{\rm Re}(s_k)>r-j+1,\ 1\le j\le r,
\]
or, in general,
\[\zeta(s_1,\ldots,s_r;b_1,\ldots,b_r)=\sum_{0\le n_1<n_2<\cdots<n_r}\prod_{j=1}^r\frac{1}{(n_j+b_j)^{s_j}}, \ \ b_j> 0,
\ \ \sum_{k=j}^r{\rm Re}(s_k)>r-j+1,\ \ 1\le j\le r;
\]
 (iii) Barnes multiple  zeta function (of one variable) (\cite{Bar1901, Bar1904, CQ1992, Mat2002a, Mat2010, Rui2000})
\[\zeta_B(s)\equiv \zeta_B(s;a_1,\ldots,a_r;b)= \sum_{n_1,n_2,\ldots,n_r=0}^{\infty}\frac{1}{(a_1n_1+a_2n_2+\cdots+a_rn_r+b)^s},\quad {\rm Re}(s)>r,
\]
where $a_1, \ldots,a_r$ and  $b$ are complex parameters satisfying
some suitable conditions (e.g., all of them have positive real
parts).

For further generalizations of the Barnes multiple  zeta function
$\zeta_B(s)$, see  \cite{Mel1902, Mel1904, Shi1977, Shi1978} for the
case of one variable, and \cite{Mat2002a}--\cite{Mat2005} for the
case of several variables. The analytically continued properties of
the above-mentioned multiple zeta functions have been investigated.
For example, the Barnes multiple zeta function $\zeta_B(s)$ extends
to a meromorphic function in $s\in\C$ and has simple poles only at
$s=1, 2,..., r$ (see \cite{Bar1904, Mat2010, Rui2000}). The main
tools of analytic continuation of the zeta functions include
integral representation, Mellin--Barnes integral formula,
Euler--Maclaurin summation, translation formula (identity),
generalized function, etc. For the recent development on the
multiple zeta functions, see the monograph \cite{Zha2016}.

In this paper, we will introduce the following three variants (a) --
(c) of the Barnes multiple zeta function, mimicking the
 Shintani double  zeta functions of the form (\cite{KTW2022, Shi1975})
\[\zeta_1(s_1,s_2)=\sum_{n_1,n_2=1}^{\infty}\frac{A(4n_1,n_2)}{n_1^{s_1}n_2^{s_2}},
\quad \zeta_2(s_1,s_2)=\sum_{n_1,n_2=1}^{\infty}\frac{A(n_1,-n_2)}{n_1^{s_1}(4n_2)^{s_2}},
\]
where $A(n_1,n_2)$ denotes the number of distinct solutions to the
quadratic congruence equation $x^2=n_2\ ({\rm mod}\ n_1).$

Let $r\ge 1$ be an integer  and for convenience let
$\alpha_1,\alpha_2,\ldots,\alpha_r,$ $a_1,a_2,\ldots, a_r, b$ be
positive real numbers.  Denote ${\bf
\alpha}=(\alpha_1,\alpha_2,\ldots,\alpha_r),$ ${\bf
a}=(a_1,a_2,\ldots, a_r)$ and
$\beta=\alpha_1+\alpha_2+\cdots+\alpha_r.$ As usual, we adopt the
notations $0!=1,$ and for real $\gamma,$
\[{\gamma\choose 0}=1,\quad {\gamma\choose n}=\frac{\gamma(\gamma-1)\cdots(\gamma-n+1)}{n!}\quad {\rm if}\ n\in\N\equiv\{1,2,\dots\}.
\]
Then, we define (see Remark 1 below for the names of moment functions)\\
(a) the multiple hyperbolic sinh-moment function
\begin{eqnarray}S_{{\bf \alpha}, {\bf a}}(s;b)=\sum_{n_1,n_2,\ldots,n_r=0}^{\infty}
\frac{{n_1+\alpha_1-1\choose n_1}\cdots{n_r+\alpha_r-1\choose
n_r}}{(a_1n_1+\cdots+a_rn_r+b)^s},\quad {\rm Re}(s)>\beta;
\end{eqnarray}
(b) the multiple hyperbolic cosh-moment function
\begin{eqnarray}C_{{\bf \alpha}, {\bf a}}(s;b)=\sum_{n_1,n_2,\ldots,n_r=0}^{\infty}
\frac{{-\alpha_1\choose n_1}\cdots{-\alpha_r\choose
n_r}}{(a_1n_1+\cdots+a_rn_r+b)^s},\quad {\rm Re}(s)>\beta;
\end{eqnarray}
(c) the multiple hyperbolic tanh-moment function
\begin{eqnarray}T_{{\bf \alpha}, {\bf a}}(s;b)=\sum_{n_1,n_2,\ldots,n_r=0}^{\infty}
\E\Big\{\frac{{-\alpha_1\choose n_1}\cdots{-\alpha_r\choose n_r}}{(a_1(n_1+V_1)+\cdots+a_r(n_r+V_r)+b)^s}\Big\},\quad {\rm Re}(s)>\beta,
\end{eqnarray}
where $\alpha_1,\ldots,\alpha_r$ are integers,
$V_j=U_{j1}+\cdots+U_{j\alpha_j},\ 1\le j\le r,$ and
$U_{j1},\ldots,U_{j\alpha_j},\ 1\le j\le r,$ are independently and
identically distributed (i.i.d.) uniform random variables on
$(0,1).$

It is seen that when all $\alpha_j=1,\ 1\le j\le r,$ the function
$S_{{\bf \alpha}, {\bf a}}(s;b)$ in (3) above reduces to the Barnes
multiple zeta function $\zeta_B(s)$ (see Corollary 1 below). We will
provide the integral representations of the above three novel
multiple zeta functions through the hyperbolic probability
distributions given by Pitman and Yor \cite{PY2003}, and investigate
their analytically continued properties including the set of
singularities. Surprisingly, unlike the previous results, the
functions $C_{{\bf \alpha}, {\bf a}}(s;b)$ and $T_{{\bf \alpha},
{\bf a}}(s;b)$ in (4) and (5) above can extend analytically to
entire functions in $s\in\C.$ The main results are stated in Section
3. Section 2 provides the needed lemmas. Finally, the proofs of the
main results are given in Section 4.
\bigskip\\
\noindent{\bf 2. Lemmas}
\medskip\\
\indent We need the following preparatory lemmas. Lemma 1 presents
the basic properties of the gamma function (see, e.g., Apostal \cite[Chapter
12]{Apo1976}). Lemma 2 is about  Liouville's extension of Dirichlet
integral theorem  (see, e.g., Sharma \cite[Chapter 2]{Sha2005}, or Gibson
\cite[page 492]{Gib1931}). The crucial Lemma 3 modifies slightly
 Cheng and Hu's \cite{CH2022} Lemma 2.1.  The
hyperbolic random variables in Lemma 4 are introduced and
investigated by Pitman and Yor \cite{PY2003}. On the other hand, for
the derivations of the hyperbolic density functions in Lemma 5, see
\cite{CH2022, PY2003}.
\medskip\\
\noindent{\bf Lemma 1} (Gamma function). {\it Let $\Gamma(s)$ be the
gamma function defined by
\[\Gamma(s)=\int_0^{\infty}x^{s-1}e^{-x}dx,\quad {\rm Re}(s)>0.
\]
Then the function $\Gamma(s)$ can extend to a meromorphic function
in the whole complex plane $\C$ and has simple poles at the points
$s=0,-1,-2,\ldots,$ with residue $(-1)^n/n!$ at $s=-n,$ where
$n=0,1,2,\ldots.$ The gamma function satisfies the equation
$\Gamma(s+1)=s\,\Gamma(s).$ It has no zeros, and hence the
reciprocal gamma function, $1/\Gamma(s)$, is an entire function in
$s\in\C$ with zeros at $s=0, -1,-2,\ldots.$ Moreover, if $x>0$ is
not an integer, then
\[\Gamma(-x)=\frac{(-1)^{[x]+1}\Gamma([x]+1-x)}{x(x-1)\cdots(x-[x])},
\]
where  $[x]$ denotes the greatest integer less than or equal to
$x.$}
\medskip\\
\noindent{\bf Lemma 2} (Liouville--Dirichlet integral). {\it Let
$m\ge 2$ be an integer and let $f$ be a Lebesgue measurable function
on $[0,1].$ Then we have
\begin{eqnarray*}
& &\idotsint\limits_{A}f(t_1+\cdots+t_{m-1})\,t_1^{s_1-1}\cdots t_{m-1}^{s_{m-1}-1}dt_1\cdots
dt_{m-1}\\
& &~~~~~~~~~=\frac{\Gamma(s_1)\cdots\Gamma(s_{m-1})}{\Gamma(s_1+\cdots+s_{m-1})}
\int_0^1f(t)\,t^{s_1+\cdots+s_{m-1}-1}dt,
\end{eqnarray*}
provided that the resulting single integral on the right exists.
Here, the parameters $s_j>0,$ $1\le j\le m-1,$ and the region of
multiple integral
\[A=\{(t_1,\ldots,t_{m-1}): t_j\ge 0,\,1\le j\le m-1,\,t_1+\cdots+t_{m-1}\le 1\}.
\]
In particular, if $f(t)=(1-t)^{s_m-1},$ where $t\in[0,1]$ and
$s_m>0,$ then
\begin{eqnarray*}
& &\idotsint\limits_{A}t_1^{s_1-1}\cdots t_{m-1}^{s_{m-1}-1}(1-t_1-\cdots-t_{m-1})^{s_m-1}dt_1\cdots
dt_{m-1}\\
& &~~~~~~~~~~=\frac{\Gamma(s_1)\cdots\Gamma(s_{m})}{\Gamma(s_1+\cdots+s_{m})},\quad s_j>0,\ 1\le j\le m.
\end{eqnarray*}}
\noindent{\bf Lemma 3} (Integral transform of a distribution). {\it
Let $X$ be a random variable with distribution $F_X$ on $\R,$
denoted $X\sim F_X,$ and let $X$ have finite moments of all orders.
For any nonzero real $c,$ define the complex-valued function
\[h_X(s;c)\equiv\int_{\R}\frac{1}{(c+ix)^s}dF_X(x)=\E\big\{(c+iX)^{-s}\big\},\quad s\in\C,
\]
where $(c+ix)^s=\exp({s{\rm Log}(c+ix)}).$ Then $h_X(s;c)$ is an
entire function in $s\in\C.$}\\
{\bf Proof.} We sketch the proof below. For each $d>0,$ define the
closed  square
\[{S}(d)=\{s=\sigma+it:\ |\sigma|\le d,\ |t|\le d,\quad \sigma, t\in \R\}.
\]
Denote the integer $n_0=[d]+1.$ Then there exists a constant
$C=C(c,d,n_0,F)$ such that
\[\int_{\R}\Big |\frac{1}{(c+ix)^s}\Big |dF_X(x)\le C\quad {\rm for\ all}\ s\in S(d).
\]
This implies that the integral transformation $h_X(s;c)$ is
absolutely convergent  and converges uniformly in each compact
subset of $\C.$ Moreover, $h_X(s,c)$ is a continuous function in
$s\in\C.$ Let $\triangle$ be any closed triangle in $\C$ and recall
that for each fixed real $x,$ $(c+ix)^{-s}$ is an entire function in
$s.$ Then we have, by applying first Fubini's theorem and then
Cauchy's theorem,
\[\int_{\partial{\triangle}}h_X(s,c)ds=\int_{\partial{\triangle}}\int_{\R}\frac{1}{(c+ix)^s}dF_X(x)ds
=\int_{\R}\int_{\partial{\triangle}}\frac{1}{(c+ix)^s}dsdF_X(x)=0,
\]
where $\partial{\triangle}$ is the boundary of $\triangle.$ Finally,
Morera's theorem completes the proof.
\medskip\\
\noindent{\bf Remark 1.} Note that for $s=-n,$ $n\in\N,$ the
function value of $h_X(s;c)$ becomes a linear combination of moments
of $X\sim F_X.$ Thus, for simplicity, we might call the integral
transformation $h_X(s;c)$ a moment function of the distribution
$F_X.$
\medskip\\
\noindent{\bf Lemma 4} (Hyperbolic random variables). {\it For
$t>0,$ let $\hat{S}_t, \hat{C}_t, \hat{T}_t$ be three random
variables with characteristic functions
\[\E(e^{i\theta\hat{S}_t})=\Big(\frac{\theta}{\sinh\theta}\Big)^t,\quad \E(e^{i\theta\hat{C}_t})=\Big(\frac{1}{\cosh\theta}\Big)^t,\quad
\E(e^{i\theta\hat{T}_t})=\Big(\frac{\tanh\theta}{\theta}\Big)^t, \quad \theta\in\R,
\]
respectively. Then the hyperbolic random variables $\hat{S}_t,
\hat{C}_t, \hat{T}_t$ are infinitely divisible and have finite
moments of all orders.}
\medskip\\
\indent
 For $t>0,$
let $\hat{S}_t\sim F_{S,t},\, \hat{C}_t\sim F_{C,t},\, \hat{T}_t\sim
F_{T,t},$ where $F_{S,t},\, F_{C,t},\, F_{T,t}$ are called
hyperbolic $\sinh$-, $\cosh$-, $\tanh$-distributions, respectively.
\medskip\\
\noindent{\bf Lemma 5} (Hyperbolic distributions). {\it For $t>0,$
$F_{S,t},\,  F_{C,t},\,  F_{T,t}$ are absolutely continuous
symmetric distributions $($about zero$)$ and have probability
density functions
\[F_{S,t}^{\prime}(x)=f_{S,t}(x)=\frac{2^t}{\pi}\int_0^1\frac{\cos(x\ln y)y^{t-1}(\ln y)^t}{(y^2-1)^t}dy,\quad x\in\R,
\]
\[F_{C,t}^{\prime}(x)=f_{C,t}(x)=\frac{2^{t-2}}{\pi}{\rm B}((t+ix)/2, (t-ix)/2),\quad x\in\R,
\]
\[F_{T,t}^{\prime}(x)=f_{T,t}(x)=\frac{1}{\pi}\int_0^1\frac{\cos(x\ln y)(y^2-1)^t}{(\ln y)^t(y^2+1)^2}\cdot\frac{1}{y}dy,\quad x\in\R,
\]
respectively, where ${\rm B}$ is the Beta function
\[{\rm B}(s_1,s_2)=\int_0^1t^{s_1-1}(1-t)^{s_2-1}dt=\frac{\Gamma(s_1)\Gamma(s_2)}{\Gamma(s_1+s_2)},\quad {\rm Re}(s_1)>0,\ \ {\rm Re}(s_2)>0.
\]
In particular, for $t=1,2,$ we have the even density functions:
\[f_{S,1}(x)=\frac{\pi}{4\cosh^2(\pi x/2)},\quad f_{S,2}(x)=\frac{\frac{\pi}{2}(\frac{\pi}{2}x\coth(\frac{\pi}{2}x)-1)}
{\sinh^2(\pi x/2)},\quad x\in\R;
\]
\[f_{C,1}(x)=\frac{1}{2\cosh(\pi x/2)},\quad f_{C,2}(x)=\frac{x}{2\sinh(\pi x/2)},\quad x\in\R;
\]
\[f_{T,1}(x)=\frac{1}{\pi}\log\coth(\frac{\pi}{4}|x|),\quad f_{T,2}(x)=\int_{|x|}^{\infty}\frac{y(y-|x|)}{2\sinh(\pi y/2)}dy, \quad x\in\R.
\]}

\noindent{\bf 3. Main results}

Throughout the section  we assume that $r\ge 1$ is an integer and
that $\alpha_1,\ldots,\alpha_r$,$ a_1, \ldots,a_r$ and $b$ are
positive real numbers satisfying $c\equiv
b-\frac12(a_1\alpha_1+a_2\alpha_2+\cdots+a_r\alpha_r)\ne 0.$ Denote
${\bf \alpha}=(\alpha_1,\ldots,\alpha_r),\,{\bf
a}=(a_1,\ldots,a_r),$ $\beta=\alpha_1+\alpha_2+\cdots+\alpha_r$ and
$c_r=(a_1^{\alpha_1}\cdots a_r^{\alpha_r})^{-1}.$ Moreover, the
notation ``$\, \stackrel{d}{=}\, $'' means equality in distribution.
Then we have the following results, which as mentioned below
generalize some of the previous results given in Cheng and Hu
\cite{CH2022}.
\medskip\\
\noindent{\bf Theorem 1.} {\it Assume that
$X_{j}\stackrel{d}{=}\frac12 \hat{S}_{\alpha_j}, j=1,2,\ldots,r,$
are independent random variables, where $\hat{S}_{\alpha_j}$ are
hyperbolic sinh random variables defined in Lemma $4$. Then the
multiple hyperbolic sinh-moment function $S_{{\bf \alpha}, {\bf
a}}(s;b)$ in $(3)$ has an integral representation of the form:
\begin{eqnarray}
S_{{\bf \alpha}, {\bf a}}(s;b)=c_r\frac{\Gamma(s-\beta)}{\Gamma(s)}\!\int_{\R}\frac{1}
{(c+iy)^{s-\beta}}dF_{Y_r}(y)=c_r\frac{\Gamma(s-\beta)}{\Gamma(s)}\E\big\{(c+iY_r)^{-s+\beta}\big\},\
{\rm Re}(s)>\beta,
\end{eqnarray}
where $Y_r=a_1X_{1}+a_2X_{2}+\cdots+a_rX_{r}.$
Moreover, the following statements hold true.\\
$(a)$ The zeta function $S_{{\bf \alpha}, {\bf a}}(s;b)$ in $(6)$
can extend
to a meromorphic function in $s\in\C$ and all the poles are simple.\\
$(b)$ If the parameter $\beta$ is a positive integer, then the
function $S_{{\bf \alpha}, {\bf a}}(s;b)$ in $(6)$ is an analytic
function in $s\in\C$ except for $s=1,2,\ldots,\beta.$  Precisely,
\[S_{{\bf \alpha}, {\bf a}}(s;b)=\frac{c_r}{(s-1)\cdots(s-\beta)}\int_{\R}\frac{1}{(c+iy)^{s-\beta}}dF_{Y_r}(y),\quad s\ne 1,2,\ldots,\beta.
\]
$(c)$  If the parameter $\beta>0$ is not an integer, then the zeta
function $S_{{\bf \alpha}, {\bf a}}(s;b)$ in $(6)$ has trivial zeros
at $s=0,-1,-2,\ldots,$ and
 we have the
following two cases for poles.\\
$(i)$ For non-negative integer $n<\beta,$  the point $\beta-n$ is a
simple pole of $S_{{\bf \alpha}, {\bf a}}(s;b)$ with residue
\[\lim_{s\to\beta-n}(s-\beta+n)S_{{\bf \alpha}, {\bf a}}(s;b)=\frac{(-1)^nc_r\E\big((c+iY_r)^{n}\big)}{n!\Gamma(\beta-n)}.
\]
$(ii)$ For positive integer $n>\beta,$  the point $\beta-n$ is a
simple pole of $S_{{\bf \alpha}, {\bf a}}(s;b)$ with residue
\[\lim_{s\to\beta-n}(s-\beta+n)S_{{\bf \alpha}, {\bf a}}(s;b)=\frac{(-1)^nc_r\E\big((c+iY_r)^{n}\big)}{n!\,d_{\beta,n}},
\]
where
\[d_{\beta,n}=\frac{(-1)^{[n-\beta]+1}\Gamma(\beta-n+1+[n-\beta])}{(n-\beta)(n-\beta-1)\cdots(n-\beta-[n-\beta])}.
\]}

\indent  If $\alpha_1=\alpha_2=\cdots=\alpha_r=1,$ then $\beta=r$
and the multiple hyperbolic sinh-moment function $S_{{\bf \alpha},
{\bf a}}(s;b)$ in (3) reduces to the Barnes multiple  zeta function
$\zeta_B(s).$ Thus, we have the following.
\medskip\\
\noindent{\bf Corollary 1.} {\it The Barnes multiple  zeta function
$\zeta_B(s)$ has the integral representation
\begin{eqnarray}\zeta_B(s)=\frac{c_r}{(s-1)\cdots(s-r)}\int_{\R}\frac{1}{(c+iy)^{s-r}}dF_{Y_r}(y),\quad {\rm Re}(s)>r,
\end{eqnarray}
where $c_r=(a_1\cdots a_r)^{-1},$ $c=b-\frac12(a_1+\cdots+a_r)\ne
0,$ $Y_r=a_1X_1+\cdots+a_rX_r$ and $X_1,\ldots,X_r$ are $r$
independent random variables, each $X_j\stackrel{d}{=}\hat{S}_1/2$
having the density function
\begin{eqnarray}f_1(x)=\frac{\pi}{2\cosh^2(\pi x)},\quad x\in\R.
\end{eqnarray}
 Moreover, $\zeta_B(s)$ in $(7)$ extends to a meromorphic function in $s\in\C$ and has simple poles
at $s=1,2,\ldots,r.$}
\medskip\\
\indent  If, further, $r=1$ and $\alpha_1=a_1=1,$ then the multiple
hyperbolic sinh-moment function $S_{{\bf \alpha}, {\bf a}}(s;b)$
reduces to the Hurwitz zeta function $\zeta(s;b).$ Thus, we have the
following.
\medskip\\
\noindent{\bf Corollary 2} (\cite{CH2022}). {\it The Hurwitz zeta
function $\zeta(s;b)$ has the integral representation
\begin{eqnarray}\zeta(s;b)=\frac{1}{s-1}\int_{\R}\frac{1}{(c+iy)^{s-1}}dF_{Y_1}(y),\quad {\rm Re}(s)> 1,
\end{eqnarray}
where $c=b-\frac12\ne 0,$ $Y_1=X_1$ and
$X_1\stackrel{d}{=}\hat{S}_1/2$ has the density function in $(8).$
 Moreover, $\zeta(s;b)$ in $(9)$ extends to a meromorphic function in $s\in\C$ and has a simple pole
at $s=1.$}
\medskip\\
\indent  If $r=2$ and $\alpha_1=\alpha_2=1, a_1=a_2=1,$ then
$\beta=2$ and the multiple hyperbolic sinh-moment function $S_{{\bf
\alpha}, {\bf a}}(s;b)$ in (3) reduces to
\begin{eqnarray}
\zeta_2(s;b)=\sum_{n_1,n_2=0}^{\infty}\frac{1}{(n_1+n_2+b)^s}=\sum_{n=1}^{\infty}\frac{n}{(n+b-1)^s},\quad {\rm Re}(s)>2.
\end{eqnarray}

\noindent{\bf Corollary 3} (\cite{CH2022}). {\it The function
$\zeta_2(s;b)$  in $(10)$ has the integral representation
\begin{eqnarray}\zeta_2(s;b)=\frac{1}{(s-1)(s-2)}\int_{\R}\frac{1}{(c+iy)^{s-2}}dF_{Y_2}(y),\quad {\rm Re}(s)> 2,
\end{eqnarray}
where $c=b-1\ne 0$ and $Y_2=X_1+X_2\stackrel{d}{=}\hat{S}_2/2$ has
the density function
\[f_2(y)=\frac{\pi(\pi y\coth(\pi y)-1)}{\sinh^2(\pi y)},\quad y\in\R.
\]
 Moreover, $\zeta_2(s;b)$ in $(11)$ extends to a meromorphic function in $s\in\C$ and has  simple poles
 at $s=1,2.$}

\noindent{\bf Theorem 2.} {\it Assume that
$X_{j}\stackrel{d}{=}\frac12 \hat{C}_{\alpha_j}, j=1,2,\ldots,r,$
are independent random variables, where $\hat{C}_{\alpha_j}$ are
hyperbolic cosh random variables defined in Lemma $4.$ Then the
multiple hyperbolic cosh-moment function $C_{{\bf \alpha}, {\bf
a}}(s;b)$ in $(4)$ has an integral representation of the form:
\begin{eqnarray}
C_{{\bf \alpha}, {\bf a}}(s;b)=2^{-\beta}\int_{\R}\frac{1}{(c+iz)^{s}}dF_{Z_r}(z)=2^{-\beta}\E\big\{(c+iZ_r)^{-s}\big\},\quad
{\rm Re}(s)>\beta,
\end{eqnarray}
where $Z_r=a_1X_{1}+a_2X_{2}+\cdots+a_rX_{r}.$ Moreover, the zeta
function $C_{{\bf \alpha}, {\bf a}}(s;b)$ in $(12)$ can extend
analytically to
 an entire function in $s\in\C.$}
\medskip\\
\indent  If $\alpha_1=\alpha_2=\cdots=\alpha_r=1,$ then $\beta=r$
and the multiple hyperbolic cosh-moment function $C_{{\bf \alpha},
{\bf a}}(s;b)$ in (4) reduces to the following zeta function
\begin{eqnarray}
\zeta_{BC}(s)=\sum_{n_1,\ldots,n_r=0}^{\infty}\frac{(-1)^{n_1+\cdots+n_r}}{(a_1n_1+\cdots+a_rn_r+b)^s},\quad {\rm Re}(s)>r
\end{eqnarray}
(compare with the zeta function $\zeta_B(s)$).  Thus, we have the
following.
\medskip\\
\noindent{\bf Corollary 4.} {\it The zeta function $\zeta_{BC}(s)$
in $(13)$ has the integral representation
\begin{eqnarray}
\zeta_{BC}(s)=2^{-r}\int_{\R}\frac{1}{(c+iz)^{s}}dF_{Z_r}(z),\quad
{\rm Re}(s)>r,
\end{eqnarray}
where $c_r=(a_1\cdots a_r)^{-1},$ $c=b-\frac12(a_1+\cdots+a_r)\ne
0,$ $Z_r=a_1X_{1}+a_2X_{2}+\cdots+a_rX_{r}$ and $X_1,\ldots,X_r$ are
$r$ independent random variables, each
$X_j\stackrel{d}{=}\hat{C}_1/2$ having the density function
\begin{eqnarray}g_1(x)=\frac{1}{\cosh(\pi x)},\quad x\in\R.
\end{eqnarray}
Moreover, the zeta function $\zeta_{BC}(s)$ in $(14)$ extends to an
entire function in $s\in\C.$}
\medskip\\
\indent  If, further, $r=1$ and $a_1=\alpha_1=1,$ then the zeta
function $\zeta_{BC}(s)$ in (13) reduces to
\begin{eqnarray}
\zeta_{C}(s)=\sum_{n=0}^{\infty}\frac{(-1)^{n}}{(n+b)^s},\quad {\rm Re}(s)>1.
\end{eqnarray}
Thus, we have the following.
\medskip\\
\noindent{\bf Corollary 5} (\cite{CH2022}). {\it The zeta function
$\zeta_{C}(s)$ in $(16)$ has the integral representation
\begin{eqnarray}
\zeta_{C}(s)=2^{-1}\int_{\R}\frac{1}{(c+iz)^{s}}dF_{Z_1}(z),\quad
{\rm Re}(s)>1,
\end{eqnarray}
where  $c=b-\frac12\ne 0$ and $Z_1=X_1$ has the density function
$g_1$ in $(15).$ Moreover, the zeta function $\zeta_{C}(s)$ in
$(17)$ extends to an entire function in $s\in\C.$}
\medskip\\
\indent If $r=2$ and $\alpha_1=\alpha_2=1, a_1=a_2=1,$ then
$\beta=2$ and the multiple hyperbolic cosh-moment function $C_{{\bf
\alpha}, {\bf a}}(s;b)$ in (4) reduces to
\begin{eqnarray}
\zeta_{C_2}(s)=\sum_{n_1,n_2=0}^{\infty}\frac{(-1)^{n_1+n_2}}{(n_1+n_2+b)^s}=\sum_{n=0}^{\infty}\frac{(-1)^n(n+1)}{(n+b)^s},\quad {\rm Re}(s)>2.
\end{eqnarray}
Then we have the following.
\medskip\\
\noindent{\bf Corollary 6} (\cite{CH2022}). {\it The zeta function
$\zeta_{C_2}(s)$ in $(18)$  has the integral representation
\begin{eqnarray}
\zeta_{C_2}(s)=\frac14\int_{\R}\frac{1}{(c+iz)^{s}}dF_{Z_2}(z),\quad
{\rm Re}(s)>2,
\end{eqnarray}
where $c=b-1\ne 0$ and $Z_2=X_1+X_2\stackrel{d}{=}\hat{C}_2/2$ has
the density function
\[g_2(z)=\frac{2z}{\sinh(\pi z)},\quad z\in\R.
\]
Moreover, the zeta function $\zeta_{C_2}(s)$ in $(19)$ extends to an
entire function in $s\in\C.$}
\medskip\\
\noindent{\bf Theorem 3.} {\it Assume that
$X_{j}\stackrel{d}{=}\frac12 \hat{T}_{\alpha_j}, j=1,2,\ldots,r,$
are independent random variables, where $\alpha_1,\ldots,\alpha_r$
are positive integers and $\hat{T}_{\alpha_j}$ are hyperbolic tanh
random variables defined in Lemma $4.$ Then the multiple hyperbolic
tanh-moment function $T_{{\bf \alpha}, {\bf a}}(s;b)$ in $(5)$ has
an integral representation of the form:
\begin{eqnarray}
T_{{\bf \alpha}, {\bf a}}(s;b)=2^{-\beta}\int_{\R}\frac{1}{(b+iw)^{s}}dF_{W_r}(w)=2^{-\beta}\E\{(b+iW_r)^{-s}\},\quad
{\rm Re}(s)>\beta,
\end{eqnarray}
where $W_r=a_1X_{1}+a_2X_{2}+\cdots+a_rX_{r}.$ Moreover, the zeta
function $T_{{\bf \alpha}, {\bf a}}(s;b)$ in $(20)$ can extend
analytically to an entire function in $s\in\C.$}
\medskip\\
\indent If $\alpha_1=\alpha_2=\cdots=\alpha_r=1,$ then $\beta=r$ and
the multiple hyperbolic tanh-moment function $T_{{\bf \alpha}, {\bf
a}}(s;b)$ in (5) reduces to
\begin{eqnarray}\zeta_{T, {\bf a}}(s;b)=\sum_{n_1,n_2,\ldots,n_r=0}^{\infty}
\E\Big\{\frac{(-1)^{n_1+\cdots+n_r}}{(a_1(n_1+V_1)+\cdots+a_r(n_r+V_r)+b)^s}\Big\},\quad {\rm Re}(s)>r,
\end{eqnarray}
where $V_j=U_{j1}+\cdots+U_{j\alpha_j} ,\ 1\le j\le r,$  and
$U_{j\ell_j},\ 1\le j\le r,\ 1\le\ell_j\le \alpha_j,$ are i.i.d.
uniform random variables on $(0,1).$
 Thus, we have the
following.
\medskip\\
\noindent{\bf Corollary 7.}  {\it The multiple hyperbolic
tanh-moment function $\zeta_{T, {\bf a}}(s;b)$ in $(21)$
 has the integral
representation \begin{eqnarray}\zeta_{T, {\bf
a}}(s;b)=\frac{1}{2^r}\int_{\R}\frac{1}{(b+iw)^{s}}dF_{W_r}(w),\quad
{\rm Re}(s)>r,
\end{eqnarray}
where  $W_r=a_1X_1+\cdots+a_rX_r$ and $X_1,\ldots,X_r$ are  i.i.d.
random variables, each $X_j\stackrel{d}{=}\hat{T}_1/2$ having the
density function
\begin{eqnarray}h_1(x)=\frac{2}{\pi}\log\coth(\frac{\pi}{2}|x|),\quad x\in\R.
\end{eqnarray}
 Moreover, $\zeta_{T, {\bf a}}(s;b)$ in $(22)$ extends to an entire function in $s\in\C.$}
\medskip\\
\indent If, further, $r=1$ and $\alpha_1=a_1=1,$ then the
 multiple hyperbolic tanh-moment
function $\zeta_{T, {\bf a}}(s;b)$ in (21) reduces to
\begin{eqnarray}
\zeta_T(s;b)=\sum_{n=0}^{\infty}\E\Big\{\frac{(-1)^n}{(n+b+U)^s}\Big\},\quad {\rm Re}(s)>1,
\end{eqnarray}
where  $U$ has the uniform distribution in $(0,1).$
\medskip\\
\noindent{\bf Corollary 8} (\cite{CH2022}). {\it The multiple
hyperbolic tanh-moment function in $(24)$
 has the integral representation
\begin{eqnarray}\zeta_{T}(s;b)=\frac{1}{2}\int_{\R}\frac{1}{(b+iw)^{s}}dF_{W_1}(w),\quad
{\rm Re}(s)>1,
\end{eqnarray}
where $W_1=X_1\stackrel{d}{=}\hat{T}_1/2$ has the density function
in $(23).$ Moreover, $\zeta_{T}(s;b)$ in $(25)$ extends to an entire
function in $s\in\C.$}
\medskip\\
\indent If $r=2$ and $\alpha_1=\alpha_2=1, a_1=a_2=1,$ then
$\beta=2$ and the multiple hyperbolic tanh-moment function
$\zeta_{T, {\bf a}}(s;b)$ in (21) reduces to
\begin{eqnarray}
\zeta_{T_2}(s;b)=\sum_{n=0}^{\infty}
\E\Big\{\frac{(-1)^{n}(n+1)}{(n+b+U_1+U_2)^s}\Big\},\quad {\rm Re}(s)>2,
\end{eqnarray}
where  $U_1,U_2$ are i.i.d. uniform random variables on $(0,1).$
\medskip\\
\noindent{\bf Corollary 9} (\cite{CH2022}). {\it The zeta function
$\zeta_{T_2}(s;b)$ in $(26)$ has the integral representation
\begin{eqnarray}\zeta_{T_2}(s;b)=\frac{1}{4}\int_{\R}\frac{1}{(b+iw)^{s}}dF_{W_2}(w),\quad {\rm Re}(s)>2,
\end{eqnarray}
where $W_2=X_1+X_2\stackrel{d}{=}\hat{T}_2/2$ has the density
function
\[h_2(w)=\int_{2|w|}^{\infty}\frac{y(y-2|w|)}{\sinh(\pi y/2)}dy, \quad w\in\R.
\]
Moreover, $\zeta_{T_2}(s;b)$ in $(27)$ extends to an entire function
in $s\in\C.$}
\medskip\\
{\bf Remark 2.} In addition to the above single variable cases, it
is also possible to consider the following cases of several
variables, where ${\bf s}=(s_1,s_2,\ldots,s_r)\in\C^r:$
\[S^*_{{\bf \alpha}, {\bf a}}({\bf s};{\bf b})=\sum_{n_1,n_2,\ldots,n_r=0}^{\infty}\prod_{j=1}^r
\frac{{n_j+\alpha_j-1\choose n_j}}{(a_1n_1+\cdots+a_rn_r+b_j)^{s_j}}, \quad {\rm Re}(s_j)>\alpha_j,\ 1\le j\le r,
\]
and
\[C^*_{{\bf \alpha}, {\bf a}}({\bf s};{\bf b})=\sum_{n_1,n_2,\ldots,n_r=0}^{\infty}\prod_{j=1}^r
\frac{{-\alpha_j\choose n_j}}{(a_1n_1+\cdots+a_rn_r+b_j)^{s_j}}, \quad {\rm Re}(s_j)>\alpha_j,\  1\le j\le r,
\]
in which the parameters $\alpha_j, a_j,\ b_j> 0,\ 1\le j\le r,$ and
${\bf \alpha}=(\alpha_1,\ldots,\alpha_r),\,{\bf
a}=(a_1,\ldots,a_r),$ ${\bf b}=(b_1,\ldots,b_r).$ In these cases,
Lemma 2 applies. The results are too complicated to be presented
here.
\bigskip\\
\noindent{\bf 4. Proofs of the main results}

It suffices to prove Theorems 1, 2 and 3, because the corollaries
can be derived easily from the theorems.

\noindent {\bf Proof of Theorem 1.} We first claim that the multiple
series on the RHS of (3) converges absolutely provided ${\rm
Re}(s)>\beta.$ To see this, note that by Stirling's formula, or more
generally, $\Gamma(x)\approx \sqrt{2\pi}x^{x-\frac12}e^{-x}$ as
$x\to\infty,$  there exist some constants $c_{\alpha_j}>0,\ 1\le
j\le r,$ such that
\[{n_j+\alpha_j-1\choose n_j}=\frac{\Gamma(n_j+\alpha_j)}{\Gamma(n_j+1)\Gamma(\alpha_j)}\approx c_{\alpha_j}n_j^{\alpha_j-1}\ \ {\rm as}\ n_j\to\infty,\ 1\le j\le r.
\]
On the other hand, for fixed $s\in\C$ with ${\rm Re}(s)>\beta,$ take
$\varepsilon=\varepsilon(s)>0$ such that ${\rm
Re}(s)=\beta+r\varepsilon,$ then
\begin{eqnarray*}(a_1n_1+\cdots+a_rn_r+b)^{{\rm
Re}(s)}=(a_1n_1+\cdots+a_rn_r+b)^{\beta+r\varepsilon}
> \prod_{j=1}^ra_j^{\alpha_j+\varepsilon}n_j^{\alpha_j+\varepsilon}.
\end{eqnarray*}
Finally, apply the fact that each series
$\sum_{n_j=1}^{\infty}{n_j^{-(1+\varepsilon)}}$ converges.

 Multiplying both sides of the
multiple hyperbolic sinh-moment function $S_{{\bf \alpha}, {\bf
a}}(s;b)$ in (3) by $\Gamma(s),$ we have
\begin{eqnarray}\Gamma(s)S_{{\bf \alpha}, {\bf a}}(s;b)=\sum_{n_1,n_2,\ldots,n_r=0}^{\infty}
& &{n_1+\alpha_1-1\choose n_1}\cdots{n_r+\alpha_r-1\choose
n_r}\nonumber\\
& &\cdot\ \ \frac{\Gamma(s)}{(a_1n_1+\cdots+a_rn_r+b)^s},\quad {\rm Re}(s)>\beta.
\end{eqnarray}
Recall the gamma identity
\begin{eqnarray}\frac{\Gamma(s)}{a^s}=\int_0^{\infty}x^{s-1}e^{-ax}dx,\quad {\rm Re}(s)>0,\quad a>0,
\end{eqnarray}
and the binomial series $(1+x)^{-r}=\sum_{k=0}^{\infty}{-r\choose
k}x^k,\quad |x|<1,$ or, equivalently,
\begin{eqnarray}
(1-x)^{-r}=\sum_{k=0}^{\infty}{-r\choose k}(-x)^k=\sum_{k=0}^{\infty}{k+r-1\choose k}x^k,\quad |x|<1.
\end{eqnarray}
Now, by (29), (30), the above claim  and Fubini's theorem, we can
rewrite (28) as follows:
\begin{eqnarray}
\Gamma(s)S_{{\bf \alpha}, {\bf a}}(s;b)&=&\sum_{n_1,n_2,\ldots,n_r=0}^{\infty}{n_1+\alpha_1-1\choose n_1}\cdots{n_r+\alpha_r-1\choose
n_r}\nonumber\\
& &~~~~~~~~~~~~~~ \cdot \int_0^{\infty}x^{s-1}e^{-(a_1n_1+\cdots+a_rn_r+b)x}dx\nonumber\\
&=&\int_0^{\infty}x^{s-1}e^{-bx}\Big\{\sum_{n_1,n_2,\ldots,n_r=0}^{\infty}{n_1+\alpha_1-1\choose n_1}\cdots{n_r+\alpha_r-1\choose
n_r}\nonumber\\
& &~~~~~~~~~~~~~~ ~~~~~~~~~~~~~\cdot e^{-(a_1n_1+\cdots+a_rn_r)x}\Big\}dx\nonumber\\
&=&\int_0^{\infty}x^{s-1}e^{-bx}\prod_{j=1}^r(1-e^{-a_jx})^{-\alpha_j}dx,
\quad {\rm Re}(s)>\beta.
\end{eqnarray}
Next, consider the product term of the integrand in (31) and write,
using Lemma 4,
\begin{eqnarray}
\prod_{j=1}^r(1-e^{-a_jx})^{-\alpha_j}&=&\prod_{j=1}^r\Big\{\frac{e^{a_jx/2}a_jx/2}{a_jx\sinh(a_jx/2)} \Big\}^{\alpha_j}
=\prod_{j=1}^r\frac{e^{a_j\alpha_jx/2}}{a_j^{\alpha_j}x^{\alpha_j}}\E\big(e^{ia_jxX_j}\big)\nonumber\\
&=&\frac{x^{-\beta}}{a_1^{\alpha_1}\cdots a_r^{\alpha_r}}\E\big(e^{(\lambda+iY_r)x}\big),\quad x>0,
\end{eqnarray}
where $\lambda=\frac12(a_1\alpha_1+\cdots+a_r\alpha_r).$

Inserting (32) in (31) yields,  by (29) and Fubini's theorem again,
\begin{eqnarray}
\Gamma(s)S_{{\bf \alpha}, {\bf a}}(s;b)&=&\int_0^{\infty}\frac{x^{s-\beta-1}}{a_1^{\alpha_1}\cdots a_r^{\alpha_r}}\E\big(e^{-(b-\lambda-iY_r)x}\big)dx\nonumber\\
&=&\frac{1}{a_1^{\alpha_1}\cdots a_r^{\alpha_r}}\int_0^{\infty}\int_{-\infty}^{\infty}x^{s-\beta-1}e^{-(b-\lambda-iy)x}dF_{Y_r}(y)dx\nonumber\\
&=&\frac{1}{a_1^{\alpha_1}\cdots a_r^{\alpha_r}}\int_{-\infty}^{\infty}\int_{0}^{\infty}x^{s-\beta-1}e^{-(b-\lambda-iy)x}dxdF_{Y_r}(y)\nonumber\\
&=&\frac{1}{a_1^{\alpha_1}\cdots a_r^{\alpha_r}}\int_{-\infty}^{\infty}\frac{\Gamma(s-\beta)}{(b-\lambda-iy)^{s-\beta}}dF_{Y_r}(y)\nonumber\\
&=&\frac{\Gamma(s-\beta)}{a_1^{\alpha_1}\cdots a_r^{\alpha_r}}\E\big((c-iY_r)^{-s+\beta}\big)\nonumber\\
&=&c_r\Gamma(s-\beta)\E\big((c+iY_r)^{-s+\beta}\big),\quad {\rm
Re}(s)>\beta,
\end{eqnarray}
where we apply the symmetric property of the hyperbolic
sinh-distribution. This proves the integral representation in (6) by
dividing both sides of (33) by $\Gamma(s).$

Recall that the integral transformation $h_{Y_r}(s)\equiv
\E\big((c+iY_r)^{-s+\beta}\big)$ is an entire function in $s\in\C$
due to Lemma 3. Therefore, the zeta function $S_{{\bf \alpha}, {\bf
a}}(s;b)$ in (6) can extend analytically to a meromorphic function
in $s\in\C$ and all the poles are simple due to Lemma 1. This proves
part (a). Parts (b) and (c) follow from part (a) and Lemma 1 again
immediately. The proof of Theorem 1 is completed.
\medskip\\
\noindent {\bf Proof of Theorem 2.} Note first that there is a
relationship between $C_{{\bf \alpha}, {\bf a}}(s;b)$ and $S_{{\bf
\alpha}, {\bf a}}(s;b);$ specifically, we can rewrite
\[C_{{\bf \alpha}, {\bf a}}(s;b)=\sum_{n_1,n_2,\ldots,n_r=0}^{\infty}(-1)^{\sum_{j=1}^rn_j}
\frac{{n_1+\alpha_1-1\choose n_1}\cdots{n_r+\alpha_r-1\choose
n_r}}{(a_1n_1+\cdots+a_rn_r+b)^s},
\]
so that the multiple series on the RHS of (4) also converges
absolutely provided ${\rm Re}(s)>\beta.$ Then proceeding along the
same lines in the proof of Theorem 1, we can prove the theorem as
follows. Multiplying both sides of the multiple hyperbolic
sinh-moment function $C_{{\bf \alpha}, {\bf a}}(s;b)$ in (4) by
$\Gamma(s),$ we have
\begin{eqnarray}\Gamma(s)C_{{\bf \alpha}, {\bf a}}(s;b)=\!\!\sum_{n_1,n_2,\ldots,n_r=0}^{\infty}
{-\alpha_1\choose n_1}\cdots{-\alpha_r\choose
n_r} \frac{\Gamma(s)}{(a_1n_1+\cdots+a_rn_r+b)^s},\ \ {\rm Re}(s)>\beta.
\end{eqnarray}
As before, using the gamma identity and the binomial series, we can
rewrite (34) as follows:
\begin{eqnarray}
\Gamma(s)C_{{\bf \alpha}, {\bf a}}(s;b)&=&\sum_{n_1,n_2,\ldots,n_r=0}^{\infty}{-\alpha_1\choose n_1}\cdots{-\alpha_r\choose
n_r} \int_0^{\infty}x^{s-1}e^{-(a_1n_1+\cdots+a_rn_r+b)x}dx\nonumber\\
&=&\int_0^{\infty}x^{s-1}e^{-bx}\Big\{\sum_{n_1,n_2,\ldots,n_r=0}^{\infty}{-\alpha_1\choose n_1}\cdots{-\alpha_r\choose
n_r} e^{-(a_1n_1+\cdots+a_rn_r)x}\Big\}dx\nonumber\\
&=&\int_0^{\infty}x^{s-1}e^{-bx}\prod_{j=1}^r(1+e^{-a_jx})^{-\alpha_j}dx,
\quad {\rm Re}(s)>\beta.
\end{eqnarray}
Next, consider the product term of the integrand in (35) and write,
using Lemma 4,
\begin{eqnarray}
\prod_{j=1}^r(1+e^{-a_jx})^{-\alpha_j}&=&\prod_{j=1}^r\Big\{\frac{e^{a_jx/2}}{2\cosh(a_jx/2)} \Big\}^{\alpha_j}
=\prod_{j=1}^r\frac{{e^{a_j\alpha_jx/2}}}{2^{\alpha_j}}\E\big(e^{ia_jxX_j}\big)\nonumber\\
&=&\frac{1}{2^{\beta}}\E\big(e^{(\lambda+iZ_r)x}\big),\quad x>0,
\end{eqnarray}
where $\lambda=\frac12(a_1\alpha_1+\cdots+a_r\alpha_r).$

Inserting (36) in (35) yields, by (29) and Fubini's theorem again,
\begin{eqnarray}
& &\Gamma(s)C_{{\bf \alpha}, {\bf a}}(s;b)=\int_0^{\infty}\frac{x^{s-1}}{2^{\beta}}\E\big(e^{-(b-\lambda-iZ_r)x}\big)dx\nonumber\\
&=&\frac{1}{2^{\beta}}\int_0^{\infty}\int_{-\infty}^{\infty}x^{s-1}e^{-(b-\lambda-iz)x}dF_{Z_r}(z)dx\nonumber\\
&=&\frac{1}{2^{\beta}}\int_{-\infty}^{\infty}\int_{0}^{\infty}x^{s-1}e^{-(b-\lambda-iz)x}dxdF_{Z_r}(z)
=\frac{1}{2^{\beta}}\int_{-\infty}^{\infty}\frac{\Gamma(s)}{(b-\lambda-iz)^{s}}dF_{Z_r}(z)\nonumber\\
&=&\frac{\Gamma(s)}{2^{\beta}}\E\big((c-iZ_r)^{-s}\big)
=\frac{1}{2^{\beta}}\Gamma(s)\E\big((c+iZ_r)^{-s}\big),\quad {\rm
Re}(s)>\beta,
\end{eqnarray}
where we apply the symmetric property of the hyperbolic
cosh-distribution. By cancelling the common term $\Gamma(s)$ in both
sides of  (37), we
 proves the integral representation in (9).

Recall again that the integral transformation $h_{Z_r}(s)\equiv
\E\big((c+iZ_r)^{-s}\big)$ is an entire function in $s\in\C$ due to
Lemma 3. Therefore, the zeta function $C_{{\bf \alpha}, {\bf
a}}(s;b)$ in (12) can extend analytically to an entire function in
$s\in\C.$ This completes the proof of Theorem 2.
\medskip\\
\noindent {\bf Proof of Theorem 3.} As in the proof of Theorem 2,
the multiple series on the RHS of (5) converges absolutely provided
${\rm Re}(s)>\beta,$ by noting that
\[a_1(n_1+v_1)+\cdots+a_r(n_r+v_r)+b>a_1n_1+\cdots+a_rn_r+b\ \ \ {\rm for}\ \  v_j>0,\ j=1,2,\dots,r.
\]
Multiplying both sides of the function $T_{{\bf \alpha}, {\bf
a}}(s;b)$ in (5) by $\Gamma(s),$ we have
\begin{eqnarray}\Gamma(s)T_{{\bf \alpha}, {\bf a}}(s;b)&=&\sum_{n_1,n_2,\ldots,n_r=0}^{\infty}
{-\alpha_1\choose n_1}\cdots{-\alpha_r\choose
n_r}
\E\Big\{\frac{\Gamma(s)}{(a_1(n_1+V_1)+\cdots+a_r(n_r+V_r)+b)^s}\Big\},\nonumber\\
& & ~~~~~~~~~~~~~~~~~~~~~~~~~~~~~~~~~~~~~~~~~~~~~~~~~~~~~\quad {\rm
Re}(s)>\beta.
\end{eqnarray}
Then, by (29) and Fubini's theorem, we can rewrite the expectation
\begin{eqnarray}
& &\E\Big\{\frac{\Gamma(s)}{(a_1(n_1+V_1)+\cdots+a_r(n_r+V_r)+b)^s}\Big\}\nonumber\\
&=&\int_0^{\infty}\cdots\int_0^{\infty}\Gamma(s)\big(a_1(n_1+v_1)+\cdots+a_r(n_r+v_r)+b\big)^{-s}dF_{V_1}(v_1)\cdots dF_{V_r}(v_r)\nonumber\\
&=&\int_0^{\infty}\cdots\int_0^{\infty}\Big(\int_0^{\infty}x^{s-1}e^{-(a_1(n_1+v_1)+\cdots+a_r(n_r+v_r)+b)x}dx\Big)dF_{V_1}(v_1)\cdots dF_{V_r}(v_r)\nonumber\\
&=&\int_0^{\infty}x^{s-1}e^{-(a_1n_1+\cdots a_rn_r+b)x}\Big(\int_0^{\infty}\cdots\int_0^{\infty}\prod_{j=1}^re^{-a_jxv_j}dF_{V_1}(v_1)\cdots dF_{V_r}(v_r)\Big)dx\nonumber\\
&=&\int_0^{\infty}x^{s-1}e^{-(a_1n_1+\cdots a_rn_r+b)x}\prod_{j=1}^r\E\big(e^{-a_jxV_j}\big)dx\nonumber\\
&=&\int_0^{\infty}x^{s-1}e^{-(a_1n_1+\cdots a_rn_r+b)x}\prod_{j=1}^r\{\frac{1-e^{-a_jx}}{a_jx}\}^{\alpha_j}dx,\quad \  {\rm
Re}(s)>\beta,
\end{eqnarray}
where we use the fact that for uniform random variable $U$ on
$(0,1),$ $\E\big(e^{-xU}\big)=\frac{1-e^{-x}}{x}, \ x>0.$

 Next, recall that
$\E\big(e^{i\theta T_{\alpha_j}}\big)=\big(\frac{\tanh
\theta}{\theta}\big)^{\alpha_j}\ {\rm for}\  \theta\in\R$ and that
\begin{eqnarray}\sum_{n_1,n_2,\ldots,n_r=0}^{\infty}
{-\alpha_1\choose n_1}\cdots{-\alpha_r\choose n_r}e^{-(a_1n_1+\cdots
a_rn_r+b)x}=e^{-bx}\prod_{j=1}^r\big(1+e^{-a_jx}\big)^{-\alpha_j}.
\end{eqnarray}
Then it follows from (38)--(40) and the definition of the random
variable $W_r$ that
\begin{eqnarray}& &\Gamma(s)T_{{\bf \alpha}, {\bf a}}(s;b)=
\int_0^{\infty}x^{s-1}e^{-bx}\prod_{j=1}^r\{\frac{1}{a_jx}\}^{\alpha_j}\prod_{j=1}^r\{\frac{1-e^{-a_jx}}{1+e^{-a_jx}}\}^{\alpha_j}dx\nonumber\\
&=&2^{-\beta}\int_0^{\infty}x^{s-1}e^{-bx}\prod_{j=1}^r\{\frac{1}{a_jx/2}\}^{\alpha_j}\prod_{j=1}^r\{\tanh(a_jx/2)\}^{\alpha_j}dx\nonumber\\
&=&2^{-\beta}\int_0^{\infty}x^{s-1}e^{-bx}\prod_{j=1}^r  \E\{e^{ia_jxX_j}\}dx
=2^{-\beta}\int_0^{\infty}x^{s-1}e^{-bx}\E(e^{ixW_r})dx\nonumber\\
&=&2^{-\beta}\int_0^{\infty}x^{s-1}e^{-bx}\int_{-\infty}^{\infty}e^{ixw}dF_{W_r}(w)dx
=2^{-\beta}\int_{-\infty}^{\infty}\int_0^{\infty}x^{s-1}e^{-(b-iw)x}dxdF_{W_r}(w)\nonumber\\
&=&2^{-\beta}\int_{-\infty}^{\infty}\frac{\Gamma(s)}{(b-iw)^s}dF_{W_r}(w)
=2^{-\beta}\Gamma(s)\E\{(b-iW_r)^{-s}\}\nonumber\\
&=&2^{-\beta}\Gamma(s)\E\{(b+iW_r)^{-s}\},\quad {\rm
Re}(s)>\beta.
\end{eqnarray}
Cancelling $\Gamma(s)$ in both sides of (41), we obtain the integral
representation in (20). Finally, Lemma 3 completes the proof of the
theorem.

{}

\end{document}